\documentclass[amsthm]{autart}    % Enable this line and disable the 
\usepackage{color}
\usepackage{subfig}%this is to put multiply plots in one figure
\usepackage{graphicx}%for multiple plots
\usepackage{amsmath,amsfonts}
\usepackage{bm}%bolds greek symbols
\begin{document}
%%-----------------------------
%%      the top matter
%%-----------------------------
\begin{frontmatter}\vspace*{0in}
\title{Control of  the Landau--Lifshitz Equation}%\thanks{...}\thanks{...}% At most 5 thanks
\author{A.N. Chow and K.A. Morris}\address{Department of Applied Mathematics, University of Waterloo, Canada\\ a29chow@uwaterloo.ca,  kmorris@uwaterloo.ca}
          
\begin{keyword}                           % Five to ten keywords,  
Asymptotic stability, Equilibrium, Lyapunov function, Nonlinear control systems, Partial differential equations            % chosen from the IFAC 
\end{keyword}                             % keyword list or with the 
                                          % help of the Automatica 
                                          % keyword wizard
\begin{abstract} 
The Landau--Lifshitz equation describes the dynamics of magnetization inside a ferromagnet. This equation is nonlinear and  has an infinite number of stable equilibria.   It is desirable to control the system from one  equilibrium to another.  A control that moves the system from an arbitrary initial state, including an  equilibrium point, to a specified equilibrium is presented.  It is proven that the second point  is an 
asymptotically stable  equilibrium of  the controlled system. The results are illustrated with some simulations. % \achg{Global results are established when the set of equilibrium points is considered.} 
%The Landau--Lifshitz equation describes the dynamics of magnetization inside a ferromagnet. This equation is nonlinear and  has an infinite number of stable equilibria. Hysteresis is present in both the Landau-Lifshitz equation and its linearization. This suggests that hysteresis is due to multiple equilibrium points, not nonlinearity.   It is desirable to control the system from one  equilibrium to another.  A control that moves the system from an arbitrary initial state, including an  equilibrium point, to a specified equilibrium is presented.  It is proven that the second point  is a globally asymptotically stable  equilibrium of  the controlled system.Hysteresis is absent in the controlled system. The results are illustrated with some simulations.
\end{abstract}

\end{frontmatter}

%%%%%%%%%%%%%%%
%%%INTRODUCTION%%%%
%%%%%%%%%%%%%%%

\section{Introduction}
The Landau--Lifshitz equation describes the magnetic behaviour within ferromagnetic structures.  This equation was originally developed to model the behaviour of domain walls, which separate magnetic regions within a ferromagnet \cite{Landau1935}. Ferromagnets are often found in memory storage devices such as hard disks, credit cards or tape recordings. Each set of data stored in a memory device is uniquely assigned to a specific stable magnetic state of the ferromagnet, and hence it is desirable to control magnetization between different stable equilibria. This is difficult due to the presence of hysteresis in the Landau--Lifshitz equation.  Due to the presence of multiple equilibria,  a particular control can lead to different magnetizations. The particular path depends on the initial state of the system and looping in the input-output map is typical \cite{Chow2013ACC,Morris2011}. %\achg{A system that exhibits hysteresis is often characterized by a looping behaviour in its input-output map; that is, given a particular input, more than one output is possible. In the case of the Landau--Lifshitz equation, this implies more than one stable equilibrium is possible for a given input \cite{Amenda_thesis,Chow2013ACC}.}

%The magnetic state of a ferromagnet can be changed by an applied magnetic field   \cite{Alouges2009, Carbou2011, Carbou2008, Carbou2009, Noh2012, Wieser2011}.  
There is now an extensive body of results on control and stabilization of linear partial differential equations  (PDE's); see for instance the books \cite{Bensoussan-book,Curtain1995,LT00_1,LT00_2} and the review paper \cite{Morris_control_handbook_rev2010}. There are much  fewer results on control and stabilization of nonlinear PDE's  and the Landau--Lifshitz equation is particularly problematic. 
% \chg{Amenda- Any  more refs? Anything on the linearization?}  
Experiments and numerical simulations demonstrating the control of domains walls in a nanowire are presented in \cite{Noh2012,Wieser2011}. In \cite{Carbou2011},  the Landau--Lifshitz equation is linearized and shown to have an unstable equilibrium. 
One control objective is to stabilize this equilibrium with the control as the average of the magnetization in one direction and zero in the other two directions. In  \cite{Carbou2008,Carbou2009}, solutions to the Landau--Lifshitz equation are shown to be arbitrarily close to domain walls given a constant control. 

Stability results are often based on linearization  \cite{CarbouLabbe2006,Carbou2006,Jizzini2011,Labbe2012}. 
However, as is well known, results with a linearized model can only  predict local stability. Furthermore, for PDEs models, stability of the linearization does not necessarily predict even local stability of the original model \cite{AM2014}.  The nonlinearity in the Landau-Lifshitz equation  is discontinuous, that is, the equation is not quasi-linear,  and there are no results that can be used to conclude that local stability of the nonlinear equation follows from the nonlinear equation \cite{alJamal2013}. Furthermore, in many applications the goal is to move the system from one equilibrium to another equilibrium and linearization is not always useful in this context.

% \achg{removed paragraph about hysteresis.}

%The dynamics in the  Landau--Lifshitz equation are quite complex and exhibit hysteretic behaviour. For example, \cite{Cowburn1999,Suess2002}  investigated via experiments the shape change  of the hysteresis loop as the structure of the nanomagnet is varied. Experiments conducted on nanowires also demonstrate hysteresis loops \cite{Noh2012}.  Numerical simulations illustrating hysteresis loops are found in  \cite{Wiele2006,Yang2011}. Hysteresis in the Landau-Lifshitz equation was further investigated in \cite{Chow2013ACC}, and was shown that the primary cause of the hysteretic behavior is the presence of multiple equilibria. This makes global stabilization using a linearization problematic.

In the next section, the uncontrolled Landau--Lifshitz equation is described.  The main result on stabilization of the Landau--Lifshitz equation is in Section~\ref{seccontrol}. In Section~\ref{secExample}, simulations illustrating the results are presented. 

%%%%%%%%%%%%%%%
\section{Landau-Lifshitz Equation}
%%%%%%%%%%%%%%%

Consider the magnetization 
$$
\mathbf m(x,t)=(m_1(x,t),m_2(x,t),m_3(x,t)),
$$  
at position $x \in [0,L]$  and time $t \geq 0 $
in a long thin ferromagnetic material of length $L>0$.  If only the exchange energy term is considered,
 the magnetization
 is modelled by the one--dimensional (uncontrolled) Landau--Lifshitz equation %\chg{ add refs mentioned by reviewer. Is LL actually derived in Guo, or is it just math?}
 \cite{Brown1963},\cite[Chapter~6]{Guo2008}%\achg{added Brown reference}
 \begin{subequations} \label{eqLLcomplete}
\begin{align}\label{eqLLGuoDing}
\frac{\partial \mathbf m}{\partial t} &= \mathbf m \times  \mathbf m_{xx}-\nu\mathbf m\times\left(\mathbf m\times\mathbf m_{xx}\right)\\
\mathbf m(x,0)&=\mathbf m_0(x)
\end{align}
where $\times$ denotes the cross product and $\nu\geq  0$ is the damping parameter, which depends on the type of ferromagnet. The Landau--Lifshitz equation sometimes includes a parameter called the gyromagnetic ratio multiplying $ \mathbf m \times  \mathbf m_{xx}$, where $\mathbf m_{xx}$ means the magnetization is differentiated with respect to $x$ twice.  This has been set to $1$ for simplicity.  For more on the damping parameter and gryomagnetic ratio, see \cite{Gilbert2004}. 
%Added "where $\mathbf m_{xx}$ means the magnetization is differentiated with respect to $x$ twice" in response to email from reader who found the paper in Arxiv.

 The Landau--Lifshitz equation is a coupled set of three nonlinear  PDEs.  It is assumed that there is no magnetic flux at the boundaries and so  Neumann boundary conditions are appropriate:
\begin{equation}\label{eqboundarycondition}
 \mathbf m_x(0,t)=  \mathbf m_x(L,t)=\mathbf 0
\end{equation}
\end{subequations}
where $\mathbf m_{x}$ means the magnetization is differentiated with respect to $x$ once. %Added "where $\mathbf m_{x}$ means the magnetization is differentiated with respect to $x$ once" in response to email from reader who found the paper in Arxiv.

Existence and uniqueness of solutions to (\ref{eqLLcomplete}) with different degrees of regularity has been shown \cite{Carbou2001,Alouges1992}.
\begin{thm}{\em \cite[Lemma~6.3.1]{Guo2008}}
If $||\mathbf m_0(x)||_{2}=1$, the solution, $\mathbf m,$ to (\ref{eqLLGuoDing}) satisfies 
 \begin{equation}\label{eqconstrain1t}
  || \mathbf m(x,t)||_{2} =1 .
 \end{equation} 
\end{thm}
The following statement is a more restrictive version of the theorem stated in \cite{Carbou2001}. 
\begin{thm}\label{thm-Carbou}
 \cite[Thm. 1.3,1.4]{Carbou2001}.
If $\mathbf m_0 \in H_2 (0,L)$, $\mathbf m_{0 ,x} (0)=\mathbf m_{0,x}(L)=\mathbf 0$ and $\| \mathbf m_0 \| =1,$ then  there exists a time $T^*>0$ and a unique solution $\mathbf m $ of (\ref{eqLLcomplete}) such that for all $T<T^*$, $\mathbf m\in \mathcal C([0,T;  H_2 (0,L) ) \cap {\mathcal L}_2 (0,L;H_3 (0,L))  .$ 
 \end{thm}

 With more general initial conditions, solutions to (\ref{eqLLcomplete}) are defined on $\mathcal L_2^3 = \mathcal L_2 ([0,L]; \mathbb R^3)$ with the usual inner--product and norm.  The notation $\|\cdot\|_{\mathcal L_2^3}$ is used for the norm.
Define the operator
\begin{equation}
f(\mathbf m)=\mathbf m \times  \mathbf m_{xx}-\nu\mathbf m\times\left(\mathbf m\times \mathbf m_{xx}\right), 
\label{defn:f}
\end{equation}
and its domain 
\begin{align}
\label{setDforfullLL}
D=\{ & \mathbf m\in \mathcal L_2^3 :  \mathbf m_x \in \mathcal L_2^3, \, \, \mathbf m_{xx} \in \mathcal L_2^3, \hspace{2em}  \nonumber \\ 
& \;  \mathbf m_x(0)=\mathbf m_x(L) = \mathbf 0    \}.  
\end{align}

%The existence and uniqueness of the solution $\mathbf m$ in (\ref{eqLLcomplete}) allows the following result to be concluded. 
\begin{thm}{\em \cite[Theorem~4.7]{Amenda_thesis}} \label{thmuncontrolledLLsemigroup}
The operator $f(\mathbf m)$ with domain  $ D$ generates a nonlinear contraction semigroup on $\mathcal L_2^3.$
\end{thm}

Ferromagnets  are magnetized to saturation \cite[Section~4.1]{Cullity2009}; that is
 $  || \mathbf m_0 (x) ||_{2} =M_s$ where  $||\cdot||_{2}$ is the Euclidean norm and $M_s$ is the magnetization saturation.  In much of the literature, $M_s$ is set to $1$; see for example, \cite[Section~6.3.1]{Guo2008},  \cite{Carbou2001,Alouges1992,Lakshmanan2011}. This convention is used here.    Physically,  this 
  means that at each point, $x$, the magnitude of $\mathbf m_0(x)$ equals the magnetization saturation.
 The initial condition $\mathbf m_0(x)$ is furthermore assumed to be real--valued.  The assumption on the initial magnetization is satisfied by the magnetization at all time.

The set of equilibrium points of (\ref{eqLLcomplete}) is  \cite[Theorem~6.1.1]{Guo2008}
\begin{align}\label{equilibriumset}
E=&\{\mathbf a=(a_1,a_2,a_3) : a_1,a_2,a_3 \mbox{ constants and }\mathbf a^\mathrm{T}\mathbf a=1 \}.
\end{align} 

\begin{thm}\label{thmlyapunov}{\em \cite[Theorem~4.11]{Amenda_thesis}}
The equilibrium set in (\ref{equilibriumset}) is asymptotically stable in the $\mathcal L_2^3$--norm.
\end{thm}

%%%%%%%%%%%%
%%Controller Design
%%%%%%%%%%%%

%This section changed completely!!
\section{Controller Design}\label{seccontrol}
%\achg{most of this section has been changed}

A control, $ \mathbf u(t)$, is introduced into the Landau-Lifshitz equation~(\ref{eqLLGuoDing}) as follows
\begin{align}
\label{eqcontrolledLL}
\frac{\partial \mathbf m}{\partial t} & = \mathbf m \times  \mathbf m_{xx}-\nu\mathbf m\times\left(\mathbf m\times \mathbf m_{xx}\right)+ \mathbf u(t)\\
 \mathbf m(x,0)&=\mathbf m_0(x).\nonumber
\end{align}
 As for the uncontrolled system, the  boundary conditions are $\mathbf m_x(0,t)=\mathbf m_x(L,t)=\mathbf 0 .$
 Equation~(\ref{eqcontrolledLL}) is the Landau--Lifshitz equation with an external magnetic field $\mathbf u.$

%asymptotically 
The goal is to choose a control  so that the system governed by the Landau--Lifshitz equation moves from an arbitrary initial condition, possibly an equilibrium point, to a specified equilibrium point $\mathbf r$.  The control needs to be chosen so that  $\mathbf r $  becomes a stable equilibrium point of the controlled system.  It can be shown that zero is an eigenvalue of the linearized uncontrolled Landau--Lifshitz equation \cite[Chapter~4.3.2]{Amenda_thesis}. For finite-dimensional linear systems,  simple proportional control  of a system with a zero eigenvalue yields asymptotic tracking of a specified  state and this motivates choosing the control
% (See Figure~\ref{figblockdiagram}.) 
\begin{equation}
\mathbf u=k(\mathbf r -\mathbf m )
\label{equ}
\end{equation}
where  $\mathbf r \in E$ is an equilibrium point  of the uncontrolled equation (\ref{eqLLcomplete}) and $k$ is a positive constant control parameter. 

The following theorem shows that for any initial condition $\mathbf m_0$ the solution to (\ref{eqcontrolledLL}) with control  $u(t) =k(\mathbf r-{\mathbf m}) $ satisfies
$\| \mathbf m(\cdot, t)\|_{\mathcal L_2^3} \leq 1.$

\begin{thm}\label{thmstrongsolutionLLcontrol}
For any $\mathbf r \in E$ define
\begin{equation}
B\mathbf m =k(\mathbf r-{\mathbf m}) \label{eqlinearAforcontrolledLL}.
\end{equation}
If $k>0$, the nonlinear operator $f+B$ with domain $D$, where $f$ and  $D$ are defined in (\ref{defn:f}), (\ref{setDforfullLL}) respectively, generates a nonlinear contraction semigroup on  ${\mathcal L_2^3}$.
\end{thm}

\begin{pf}
It will be shown that (i) $f+B$ is dissipative and  (ii) the range of $\mathrm{I}-\alpha (f+B )$ for all $\alpha>0$ is ${\mathcal L_2^3}$. 

%\chg{Shortened}
For any $\mathbf m, \mathbf y \in D$,
\begin{align*}
&\langle f(\mathbf m)+B\mathbf m -  (f(\mathbf y) +B\mathbf y),\mathbf m -\mathbf y \rangle_{\mathcal L_2^3}  \\
%&= \langle f(\mathbf m)-f(\mathbf y) +B\mathbf m-B\mathbf y,\mathbf m -\mathbf y \rangle_{\mathcal L_2^3} \\
&=\langle f(\mathbf m)-f(\mathbf y), \mathbf m -\mathbf y \rangle_{\mathcal L_2^3}  +\langle B\mathbf m-B\mathbf y,\mathbf m -\mathbf y \rangle_{\mathcal L_2^3}.
\end{align*}
Since $f$ generates a nonlinear contraction semigroup (Theorem~\ref{thmuncontrolledLLsemigroup}), $f$ is dissipative \cite[Proposition~2.98]{Luo1999} and hence
\[
\langle f(\mathbf m)-f(\mathbf y), \mathbf m -\mathbf y \rangle_{\mathcal L_2^3}\leq 0.
\]
It follows that 
\begin{align*}
&\langle f(\mathbf m)+B\mathbf m -  (f(\mathbf y) +B\mathbf y),\mathbf m -\mathbf y \rangle_{\mathcal L_2^3}  \\
&\leq \langle B\mathbf m-B\mathbf y,\mathbf m -\mathbf y \rangle_{\mathcal L_2^3}\\&=\langle -k\mathbf m+k\mathbf y,\mathbf m -\mathbf y \rangle_{\mathcal L_2^3}\\
&=-k\langle \mathbf m-\mathbf y,\mathbf m -\mathbf y \rangle_{\mathcal L_2^3}\\
&\leq 0
\end{align*}
and hence $f+B$ is dissipative.

%To show the range of $\mathrm{I}-\alpha (B+f )$ is the entire space ${\mathcal L_2^3}$, consider the following. 
Since $f$ generates a nonlinear contraction semigroup (Theorem~\ref{thmuncontrolledLLsemigroup}), then it is $m$-dissipative and hence, $\mathrm{ran}(\mathrm{I}-\hat \alpha f)=\mathcal L_2^3$ for any $\hat\alpha >0$  \cite[Lemma~2.1]{Kato1967}.  
This means  that for any $\mathbf y_2 \in \mathcal L_2^3$ there exists $\mathbf m \in D$ such that $\mathbf{m}-\hat{\alpha}f(\mathbf m)=\mathbf y_2. $ 
Choose any $\mathbf y_1 \in \mathcal L_2^3$, $\alpha >0$ and define 
\[
\mathbf y_2 = \frac{\mathbf y_1}{1+\alpha k}+\frac{\alpha k\mathbf r }{1+\alpha b k}
\]
and
\[
\hat \alpha=\frac{\alpha  }{1+\alpha k}.
\]
There exists $\mathbf m \in D$ such that
\[
\mathbf{m}-\frac{\alpha }{1+\alpha k}f(\mathbf m)=\mathbf y_2=  \frac{\mathbf y_1}{1+\alpha k}+\frac{\alpha k\mathbf r }{1+\alpha b k}.
\]
Solving for $\mathbf y_1$ leads to
\[
\mathbf y_1 = \mathbf m-\alpha (k(\mathbf r-\mathbf m)+ f(\mathbf m)).
\]
Thus, for any $\mathbf y_1 \in \mathcal L_2^3$, there exists $\mathbf m \in D$ such that $
\mathbf y_1 = (\mathrm I-\alpha (B+f))\mathbf m$ and hence $\mathrm{ran}\left(\mathrm I-\alpha (B+f)\right)=\mathcal L_2^3$ for some $\alpha >0$. It follows that the range is $\mathcal L_2^3$  for all  $\alpha>0$ \cite[Lemma~2.1]{Kato1967}. 

Thus, since $B+f$ is dissipative and the range of $(\mathrm I-\alpha (B+f))$ is $\mathcal L_2^3$, then $B+f$ generates a nonlinear contraction semigroup \cite[Proposition~2.114]{Luo1999}.
 \qed
 \end{pf}

%The following lemmas are needed in the proof of Theorem~\ref{thmr0isasymstable}. 
The results in Lemma~\ref{thmderivativemcrossmprime} and~\ref{thmderivativemcrossmprimedotproduct} are a consequence of the product rule.

\begin{lem}\label{thmderivativemcrossmprime}
For $\mathbf m\in \mathcal L_2^3$, the derivative of $\mathbf g=\mathbf m \times \mathbf m_x$ is $\mathbf g_x=\mathbf m \times \mathbf m_{xx}$.
\end{lem}

\begin{lem}\label{thmderivativemcrossmprimedotproduct}
For $\mathbf m \in \mathcal L_2^3$, the derivative of $f=\left(\mathbf m \times \mathbf m_x\right)^{\mathrm T}\left(\mathbf m \times \mathbf m_x\right)$ is $f_x=2\left(\mathbf m \times \mathbf m_{x}\right)^{\mathrm T}\left(\mathbf m \times \mathbf m_{xx}\right).$
%\[
%f_x=2\left(\mathbf m \times \mathbf m_{x}\right)^{\mathrm T}\left(\mathbf m \times \mathbf m_{xx}\right).
%\]
\end{lem}

\begin{lem}
\label{lemmazerointegral}
For $\mathbf m \in \mathcal L_2^3$ satisfying (\ref{eqboundarycondition}), 
\[
\int_0^L (\mathbf m-\mathbf r)^{\mathrm T}(\mathbf m \times \mathbf m_{xx})dx=0.
\]
\end{lem}
\begin{pf}
Integrating by parts, and applying  Lemma~\ref{thmderivativemcrossmprime} and the boundary conditions (\ref{eqboundarycondition}) implies that 
\[
\int_0^L (\mathbf m-\mathbf r)^{\mathrm T}(\mathbf m \times \mathbf m_{xx})dx = -\int_0^L \mathbf m_x^{\mathrm T}(\mathbf m \times \mathbf m_{x})dx.
\]
From properties of cross products,
$
\mathbf m_x^{\mathrm T}(\mathbf m \times \mathbf m_{x}) =  \mathbf m^{\mathrm T}(\mathbf m_x \times \mathbf m_{x})=0,
$
and hence the integral is zero. \qed
\end{pf}
\begin{lem}\label{lemmaPoincareInequalityforCrossProducts}
For $\mathbf m \in \mathcal L_2^3$ satisfying (\ref{eqboundarycondition}), %then
\[
|| \mathbf m \times \mathbf m_x ||_{\mathcal L_2^3} \leq 4L^2|| \mathbf m \times \mathbf m_{xx} ||_{\mathcal L_2^3}
\]
\end{lem}
\begin{pf}
Integrating by parts, using Lemma~\ref{thmderivativemcrossmprimedotproduct} and the boundary conditions (\ref{eqboundarycondition}) leads to
\[
|| \mathbf m \times \mathbf m_x ||_{\mathcal L_2^3}^2 =-\int_0^L2\left(\mathbf m \times \mathbf m_{x}\right)^{\mathrm T}\left(\mathbf m \times \mathbf m_{xx}\right)xdx.
\]
 It follows from Young's inequality that
\begin{align*}
|| \mathbf m \times \mathbf m_x ||_{\mathcal L_2^3}^2 %&=-\int_0^L2\left(\mathbf m \times \mathbf m_{x}\right)^{\mathrm T}\left(\mathbf m \times \mathbf m_{xx}\right)xdx\\
 &\leq \frac{1}{2}\int_0^L\left(\mathbf m \times \mathbf m_{x}\right)^{\mathrm T}\left(\mathbf m \times \mathbf m_{x}\right)dx\\
 &+\int_0^L2\left(\mathbf m \times \mathbf m_{xx}\right)^{\mathrm T}\left(\mathbf m \times \mathbf m_{xx}\right)x^2dx .
 \end{align*}
 Since $x \in [0,L],$ 
 \begin{align*} 
 &|| \mathbf m \times \mathbf m_x ||_{\mathcal L_2^3}^2
 &\leq \frac{1}{2}||\mathbf m \times \mathbf m_{x}||_{\mathcal L_2^3}^{2}+2L^2||\mathbf m \times \mathbf m_{xx}||_{\mathcal L_2^3}^2.
\end{align*}
Rearranging gives the desired inequality. \qed
\end{pf}

\begin{thm}\label{thmr0isexpstable}
Let $\mathbf r$ be an equilibrium point of (\ref{eqcontrolledLL}). For any  constant $k$ such that ${k> 8\nu L^4}$, $\mathbf r$ is  a globally exponentially stable equilibrium point of (\ref{eqcontrolledLL}) in the $H_1^3$--norm. That is, for any initial condition on $H_1^3$, $\mathbf m$ decreases exponentially in the $H_1$-norm to $\mathbf r.$
% where $||\mathbf m ||^2_{H_1}=||\mathbf m||_{\mathcal L_2^3}^2+||\mathbf m_x||_{\mathcal L_2^3}^2.$
\end{thm}

\begin{pf} 
If the initial condition $m_0 (x) \in H_1^3,$ then solutions to (\ref{eqcontrolledLL}) are in $H_1^3$ \cite{Alouges1992}. The Lyapunov candidate 
\[
V(\mathbf m)=\frac{1}{2}\left| \left| \mathbf m-\mathbf r\right|\right|_{\mathcal L_2^3}^2+\frac{1}{2}\left| \left|  \mathbf m_x\right|\right|_{\mathcal L_2^3}^2
\]
is clearly positive definite for all $\mathbf m\in D$. Furthermore, $V(\mathbf m)=0$ only when $\mathbf m=\mathbf r$. Taking the derivative of $V$ along trajectories of the controlled equation (\ref{eqcontrolledLL})
\begin{align*}
\frac{dV}{dt}&=\int_0^L(\mathbf m -\mathbf r)^{\mathrm T}\dot{{\mathbf m}} dx+\int_0^L\mathbf m_x^{\mathrm T} \dot{{\mathbf m}}_xdx \\
&=\int_0^L(\mathbf m -\mathbf r)^{\mathrm T}\dot{\mathbf m} dx-\int_0^L\mathbf m_{xx}^{\mathrm T} \dot{\mathbf m}dx
\end{align*}
where the dot notation means differentiation with respect to $t$. Substituting in  (\ref{eqcontrolledLL}) to eliminate $\dot{\mathbf m}$,
\begin{align*}
\frac{dV}{dt}
&= \int_0^L(\mathbf m -\mathbf r)^{\mathrm T}  \left(\mathbf m \times  \mathbf m_{xx}\right)dx\\
&-\nu\int_0^L (\mathbf m -\mathbf r)^{\mathrm T}\left(\mathbf m\times\left(\mathbf m\times \mathbf m_{xx}\right)\right)dx \\
&+k\int_0^L(\mathbf m -\mathbf r)^{\mathrm T} (\mathbf r -\mathbf m)dx \\
&-\int_0^L\mathbf m_{xx}^{\mathrm T}  \left(\mathbf m \times  \mathbf m_{xx}\right)dx\\
&+\nu\int_0^L \mathbf m_{xx}^{\mathrm T}\left(\mathbf m\times\left(\mathbf m\times \mathbf m_{xx}\right)\right)dx\\
& -k\int_0^L\mathbf m_{xx}^{\mathrm T} (\mathbf r -\mathbf m)dx.
\end{align*}
From Lemma~\ref{lemmazerointegral}, the first integral is zero. Furthermore, from properties of cross products, %it follows that
\[
\mathbf m_{xx}^{\mathrm T}  \left(\mathbf m \times  \mathbf m_{xx}\right)=\mathbf m^{\mathrm T}  \left(\mathbf m_{xx} \times  \mathbf m_{xx}\right)=0,
\] 
and hence
\[
\int_0^L\mathbf m_{xx}^{\mathrm T}  \left(\mathbf m \times  \mathbf m_{xx}\right)dx=0.
\]
It follows that 
\begin{align}
\frac{dV}{dt}=&-\nu\int_0^L (\mathbf m -\mathbf r)^{\mathrm T}\left(\mathbf m\times\left(\mathbf m\times \mathbf m_{xx}\right)\right)dx  \nonumber\\
& -k||\mathbf m -\mathbf r||_{\mathcal L_2^3}^2-\nu||\mathbf m\times \mathbf m_{xx}||_{\mathcal L_2^3}^2 \nonumber \\
& -k\int_0^L\mathbf m_{xx}^{\mathrm T} (\mathbf r -\mathbf m)dx. \label{dVdt} 
\end{align}
Applying integration by parts to the last integral, 
\begin{equation} \int_0^L\mathbf m_{xx}^{\mathrm T} (\mathbf r -\mathbf m)dx=||\mathbf m_{x}||_{\mathcal L_2^3}^2. \label{dVdt3} \end{equation}
%\begin{align}
%\frac{dV}{dt}=&-\nu\int_0^L (\mathbf m -\mathbf r)^{\mathrm T}\left(\mathbf m\times\left(\mathbf m\times \mathbf m_{xx}\right)\right)dx\nonumber\\
%& -bk||\mathbf m -\mathbf r||_{\mathcal L_2^3}^2-\nu||\mathbf m\times \mathbf m_{xx}||_{\mathcal L_2^3}^2 -bk||\mathbf m_{x}||_{\mathcal L_2^3}^2 \nonumber\\
%=&-\nu\int_0^L\left( (\mathbf m -\mathbf r)\times \mathbf m\right)^{\mathrm T}\left(\mathbf m\times \mathbf m_{xx}\right)dx\nonumber\\
%& -bk||\mathbf m -\mathbf r||_{\mathcal L_2^3}^2-\nu||\mathbf m\times \mathbf m_{xx}||_{\mathcal L_2^3}^2 -bk||\mathbf m_{x}||_{\mathcal L_2^3}^2 \nonumber\\
%=&\nu\int_0^L\left( \mathbf r\times \mathbf m\right)^{\mathrm T}\left(\mathbf m\times \mathbf m_{xx}\right)dx \nonumber\\
%&-bk||\mathbf m -\mathbf r||_{\mathcal L_2^3}^2-\nu||\mathbf m\times \mathbf m_{xx}||_{\mathcal L_2^3}^2 -bk||\mathbf m_{x}||_{\mathcal L_2^3}^2. \label{eqderivativeofVlastintegral}
%\end{align}
The first integral in $\frac{dV}{dt}$ can be rewritten 
\begin{align*}
-\nu\int_0^L  & (\mathbf m -\mathbf r)^{\mathrm T}\left(\mathbf m\times\left(\mathbf m\times \mathbf m_{xx}\right)\right)dx \\ 
=&-\nu\int_0^L\left( (\mathbf m -\mathbf r)\times \mathbf m\right)^{\mathrm T}\left(\mathbf m\times \mathbf m_{xx}\right)dx\\
=&\nu \int_0^L\left( \mathbf r\times \mathbf m\right)^{\mathrm T}\left(\mathbf m\times \mathbf m_{xx}\right)dx
\end{align*}
Applying integration by parts to this integral,  using Lemma~\ref{thmderivativemcrossmprime} and boundary  conditions  (\ref{eqboundarycondition}) leads to 
%Applying integration by parts with Lemma~\ref{thmderivativemcrossmprime} to the integral implies
%\begin{align*}
%&\int_0^L\left( \mathbf r\times \mathbf m\right)^{\mathrm T}\left(\mathbf m\times \mathbf m_{xx}\right)dx\\
%&=\left[\left(\mathbf r \times \mathbf m\right)^{\mathrm T}\left(\mathbf m \times \mathbf m_x\right)\right]_0^L-\int_0^L\left( \mathbf r\times \mathbf m_x\right)^{\mathrm T}\left(\mathbf m\times \mathbf m_{x}\right)dx
%\end{align*}
%and substituting in the boundary conditions in (\ref{eqboundarycondition}) leads to 
\begin{align*}
\nu\int_0^L\left( \mathbf r\times \mathbf m\right)^{\mathrm T}\left(\mathbf m\times \mathbf m_{xx}\right)dx&=-\nu\langle  \mathbf r\times \mathbf m_x, \mathbf m\times \mathbf m_{x}\rangle_{\mathcal L_2^3}.
\end{align*}
Then from Cauchy-Schwarz, Lemma~\ref{lemmaPoincareInequalityforCrossProducts}, Young's Inequality and Lemma \ref{thmderivativemcrossmprimedotproduct},
\begin{align*}
\nu\int_0^L\left( \mathbf r\times \mathbf m\right)^{\mathrm T}\left(\mathbf m\times \mathbf m_{xx}\right)dx&\leq\nu|| \mathbf r\times \mathbf m_x||_{\mathcal L_2^3} ||\mathbf m\times \mathbf m_{x}||_{\mathcal L_2^3}
\\
&\leq 4\nu L^2  || \mathbf r\times \mathbf m_x||_{\mathcal L_2^3} ||\mathbf m\times \mathbf m_{xx}||_{\mathcal L_2^3} \\
%\end{align*}
%It follows from Young's Inequality that
%\begin{align*}
%&\int_0^L\left( \mathbf r\times \mathbf m\right)^{\mathrm T}\left(\mathbf m\times \mathbf m_{xx}\right)dx\\
%&\leq 8L^4 || \mathbf r\times \mathbf m_x||_{\mathcal L_2^3} ^2+\frac{1}{2}||\mathbf m\times \mathbf m_{xx}||_{\mathcal L_2^3}^2\\
&\leq {8\nu L^4} ||  \mathbf m_x||_{\mathcal L_2^3} ^2+\frac{\nu}{2} ||\mathbf m\times \mathbf m_{xx}||_{\mathcal L_2^3}^2 .
\end{align*}
%where the latter inequality follows from Lemma \ref{thmderivativemcrossmprimedotproduct}.
Substituting this inequality and (\ref{dVdt3}) into (\ref{dVdt}) leads to
\begin{align*}%\label{eqderivativeofVLL}
\frac{dV}{dt}\leq& -\left(k- {8\nu L^4}\right)||\mathbf m_{x}||_{\mathcal L_2^3}^2\nonumber\\
&-\frac{\nu}{2}||\mathbf m\times \mathbf m_{xx}||_{\mathcal L_2^3}^2-k||\mathbf m -\mathbf r||_{\mathcal L_2^3}^2\\
&\leq -\left(k-{8\nu L^4}\right)||\mathbf m_{x}||_{\mathcal L_2^3}^2-k||\mathbf m -\mathbf r||_{\mathcal L_2^3}^2\\
&\leq -\left(k-{8\nu L^4}\right)\left(||\mathbf m_{x}||_{\mathcal L_2^3}^2+||\mathbf m -\mathbf r||_{\mathcal L_2^3}^2\right)\\
&=-2\left(k-{8\nu L^4}\right)V.
\end{align*}
Integrating with respect to time and noting that $\mathbf r$ does not depend on $x$
\begin{align*}
&||\mathbf m_{x}||_{\mathcal L_2^3}^2+||\mathbf m -\mathbf r||_{\mathcal L_2^3}^2\\
& \leq e^{-2\left(k-{8\nu L^4}\right)t}\left(||\mathbf {m}_{x}(x,0)||_{\mathcal L_2^3}^2+||\mathbf m(x,0) -\mathbf r||_{\mathcal L_2^3}^2 \right)\\
%\end{align*}
%Since $\mathbf r$ does not depend on $x$,  %which leads to
%\begin{align*}
%&||\left(\mathbf m-\mathbf r\right)_{x}||_{\mathcal L_2^3}^2+||\mathbf m -\mathbf r||_{\mathcal L_2^3}^2 \\
&= e^{-2\left(k-{8\nu L^4}\right)t}\left(||\left(\mathbf m(x,0)-\mathbf r\right)_{x}||_{\mathcal L_2^3}^2+||\mathbf m(x,0) -\mathbf r||_{\mathcal L_2^3}^2 \right).
\end{align*}
Therefore, 
\[
||\mathbf m-\mathbf r||_{H_1}^2\leq e^{-2\left(k-{8\nu L^4}\right)t}||\mathbf m(x,0) -\mathbf r||_{H_1}^2
\]
and since $k>{8\nu L^4}$, $\mathbf r$ is an exponentially stable equilibrium point of (\ref{eqcontrolledLL}). \qed
\end{pf}
%The results in Theorems~\ref{thmr0isasymstable} and~\ref{thmr0isexpstable} provide global asymptotic stability in the $\mathcal L_2^3$-norm and  global %exponential stability in the $H_1$-norm of the given equilibrium point $\mathbf r.$

%%%%%%%%%%%%%%%
%LINEARIZATION SECTION
%%%%%%%%%%%%%%%%
A natural question is whether  $\mathbf r$ is exponentially stable in the $\mathcal L_2^3$-norm. 
Analysis of the linear Landau--Lifshitz equation provides insight to this question. The linearized controlled Landau--Lifshitz equation is 
%\begin{align*}
%\frac{\partial \mathbf z}{\partial t} &=\nu \mathbf z_{xx}+\mathbf a \times {\mathbf z}_{xx} +b\mathbf u(t)\\
%\mathbf y&=\mathbf z\\
%\mathbf e&=\mathbf r-\mathbf y \\
%\mathbf u(t)&=k\mathbf e 
%\end{align*}
%where $\mathbf r$ is an equilibrium point described in Theorem~\ref{thmlinLLstable}. 
\begin{align}\label{eqcontrolledlinearLL}
\frac{\partial \mathbf z}{\partial t} & =\nu \mathbf  z_{xx}+\mathbf a \times  {\mathbf z}_{xx} +k \left(\mathbf r-\mathbf z\right), \qquad \mathbf z(0)=\mathbf z_0
\end{align}
with the same  boundary conditions $\mathbf z_{x}(0)=\mathbf z_{x}(L)=\mathbf 0 .$ 
 Since the uncontrolled linear Landau--Lifshitz equation generates a linear semigroup and $k \left(\mathbf r-\mathbf z\right)$ is a bounded linear (affine) operator, then the operator in (\ref{eqcontrolledlinearLL}) generates a semigroup \cite[Theorem~3.2.1]{Curtain1995}.  Substituting $\mathbf z=\mathbf{r}$ into (\ref{eqcontrolledlinearLL}) leads to ${\partial \mathbf z}/{\partial t}=\mathbf 0$ and hence $\mathbf{r}$ is a stable  equilibrium point of (\ref{eqcontrolledlinearLL}). 
  
\begin{thm}\label{thmlinearr0isasymstable}
Let $\mathbf r \in E .$  For any positive constant $k$, $\mathbf r$ is an exponentially stable equilibrium of the linearized system (\ref{eqcontrolledlinearLL})  in
 $\mathcal L_2^3$--norm. 
\end{thm}

\begin{pf}
For $\mathbf z \in D(A)$, where $D(A)=D$ as in equation~(\ref{setDforfullLL}), consider the Lyapunov candidate
\begin{equation*}
V(\mathbf z)=\frac{1}{2}\left| \left|\mathbf z-\mathbf r\right|\right|_{\mathcal L_2^3}^2.
\end{equation*}
It is clear that $V\geq 0$ for all $\mathbf z\in D(A)$ and furthermore, $V(\mathbf z)=0$ only when $\mathbf z=\mathbf r$.  Therefore, $V(\mathbf z)>0$ for all $\mathbf z \in D(A)\backslash  \{\mathbf r\}$.

Taking the derivative of $V(\mathbf z)$ implies
\begin{align*}
\frac{dV}{dt}&=\int_0^L(\mathbf z-\mathbf r)^{\mathrm T}\dot{ {\mathbf z}} dx.
\end{align*}
Substituting in (\ref{eqcontrolledlinearLL}) yields 
\begin{align*}
\frac{dV}{dt}&=\nu\int_0^L(\mathbf z-\mathbf r)^{\mathrm T} {\mathbf z}_{xx}dx +\int_0^L(\mathbf r-\mathbf z)^{\mathrm T}\left({\mathbf a \times {\mathbf  z_{xx}}}\right)dx \\
&+k\int_0^L(\mathbf z-\mathbf r)^{\mathrm T} ({\mathbf r-\mathbf z}) dx.
\end{align*}
By Lemma~\ref{lemmazerointegral}, the middle term is zero. Using integration by parts,  the first term becomes
\[
-\nu\int_0^L\mathbf z_x^{\mathrm T} {\mathbf z}_{x}dx.
\]
It follows that 
\begin{align*}
\frac{dV}{dt}&=-\nu||\mathbf z_x||_{\mathcal L_2^3}^2 -k||\mathbf z-\mathbf r||_{\mathcal L_2^3}^2
\end{align*}
and since $\nu\geq0$, 
\begin{align*}%\label{eqVforLLexp}
\frac{dV}{dt}&\leq -k||\mathbf z-\mathbf r||_{\mathcal L_2^3}^2=-2kV.
\end{align*}
Solving yields
\[
||\mathbf z-\mathbf r||_{\mathcal L_2^3}^2 \leq  e^{-2kt}||\mathbf z_0-\mathbf r||_{\mathcal L_2^3}^2.
\]
For $k>0$ the equilibrium point, $\mathbf r,$ of (\ref{eqcontrolledlinearLL}) is locally exponentially stable. This is true for any initial condition and hence global stability is obtained.\qed
\end{pf}
Theorem~\ref{thmlinearr0isasymstable} suggests that the equilibrium point in the controlled nonlinear Landau--Lifshitz equation (\ref{eqcontrolledLL}) is exponentially stable in the $\mathcal L_2^3$--norm. However, since the nonlinearity in the Landau-Lifshitz equation is unbounded, stability of the linear equation does not necessarily reflect stability of the original nonlinear equation; see  \cite{AM2014,alJamal2013}. 

%%%%%%%%%%%%%
%PHYSICAL CONTROL
%%%%%%%%%%%%%%
In equation~(\ref{eqcontrolledLL}), the control is affine. However, in current applications, the control enters as an applied magnetic field  \cite{Carbou2011,Carbou2008,Carbou2009,CarbouLabbe2006,Carbou2006}. More precisely,
\begin{align}
\frac{\partial \mathbf m}{\partial t}  =& \mathbf m \times \left( \mathbf m_{xx}+\mathbf u\right)-\nu\mathbf m\times\left(\mathbf m\times \left(\mathbf m_{xx}+\mathbf u\right)\right)\nonumber\\
=&\mathbf m \times  \mathbf m_{xx}-\nu\mathbf m\times\left(\mathbf m\times \mathbf m_{xx}\right) \nonumber\\
&+ \mathbf m \times  \mathbf u-\nu\mathbf m\times\left(\mathbf m\times \mathbf u\right)\label{eqcontrolledLLphysical}
\end{align}
where $\mathbf u=k(\mathbf r -\mathbf m)$ as before and $\mathbf r \in E$ is an equilibrium point of (\ref{eqcontrolledLLphysical}). Equation~(\ref{eqcontrolledLLphysical}) is the Landau-Lifshitz equation with a nonlinear control. Its existence and uniqueness results can be found in \cite[Thm. 1.1,1.2]{Carbou2001} and is similar to Theorem~\ref{thm-Carbou}.

As for the uncontrolled equation, since 
\begin{align*}
\frac{1}{2}\frac{\partial || \mathbf m(x,t)||_{2}}{\partial t}=& \mathbf m^{\mathrm T}\frac{\partial \mathbf m}{\partial t}\\
 =& \mathbf m^{\mathrm T}(\mathbf m \times  \mathbf m_{xx}-\nu\mathbf m\times\left(\mathbf m\times \mathbf m_{xx}\right) \\
&+ \mathbf m \times  \mathbf u-\nu\mathbf m\times\left(\mathbf m\times \mathbf u\right)) =0,
\end{align*}
this implies $ || \mathbf m||_{2} =c$, where $c$ is a constant. The convention is to take $c=1$. It follows that any equilibrium point is trivially stable in the $\mathcal L_2$-norm. 

%\chg{Amenda, add something about well-posedness of \ref{eqcontrolledLLphysical}. That is, the equation has a solution in some sense. Are they in $\mathcal L_2$ or $H_1$?  }

 \begin{thm}\label{thmrisstablephysical}
For any $\mathbf r \in E$ and any positive constant $k$, $\mathbf r$ is a locally stable equilibrium point of (\ref{eqcontrolledLLphysical})  in the 
%$\mathcal L_2^3$--norm.
$H_1$-norm.  That is, for any initial condition $m_0 (x) \in D,$ where $D$ is defined in (\ref{setDforfullLL}), the $H_1$-norm of the error $\mathbf m - \mathbf r$ does not increase. 
%\chg{Amenda, changed this to $H_1$ .Why would it be $L_2$ norm if the Lyapunov function is the $H_1$ norm?}
\end{thm}%%asymptotically 

\begin{pf}
%\chg{shortened}
Let $B(\mathbf r,p)=\{\mathbf m \in \mathcal L_2^3: ||\mathbf m -\mathbf r||_{ \mathcal L_2^3}<p \} \subset D$ for some constant $0<p<2$.  Note that since $p<2$, then  $-\mathbf r \notin B(\mathbf r,p)$. For any $\mathbf m \in B(\mathbf r,p)$, consider the $H_1$-norm of the error 
\[
V(\mathbf m)=k\left| \left| \mathbf m-\mathbf r\right|\right|_{\mathcal L_2^3}^2+\left| \left|  \mathbf m_x\right|\right|_{\mathcal L_2^3}^2 . 
\]
%is clearly positive definite for $k>0$. 
Taking the derivative of $V$,
\begin{align}
\frac{dV}{dt}&=\int_0^Lk(\mathbf m -\mathbf r)^{\mathrm T}\dot{{\mathbf m}} dx+\int_0^L\mathbf m_x^{\mathrm T} \dot{{\mathbf m}}_xdx \nonumber\\
&=\int_0^Lk(\mathbf m -\mathbf r)^{\mathrm T}\dot{\mathbf m} dx-\int_0^L\mathbf m_{xx}^{\mathrm T} \dot{\mathbf m}dx\nonumber\\
&=\int_0^L\left(k(\mathbf m -\mathbf r)^{\mathrm T}\dot{\mathbf m} -\mathbf m_{xx}^{\mathrm T} \dot{\mathbf m}\right)dx.\label{eqLyapunovFunction}
\end{align}

Let $\mathbf h = \mathbf m - \mathbf r$, then the integrand becomes
\begin{equation}\label{eqintegrand}
k\mathbf h^{\mathrm T}\dot{\mathbf m} - \mathbf m_{xx}^{\mathrm T} \dot{\mathbf m}
\end{equation}
and equation~(\ref{eqcontrolledLLphysical}) becomes
\[
\dot{ \mathbf m}  = \mathbf m \times \left( \mathbf m_{xx}-k\mathbf h\right)-\nu\mathbf m\times\left(\mathbf m\times \left(\mathbf m_{xx}-k\mathbf h\right)\right).
\]
It follows that %Taking the derivative of $V$ along trajectories of the controlled equation (\ref{eqcontrolledLL}) leads to  
\begin{align}
\mathbf h^{\mathrm T}\dot{\mathbf m} =& \mathbf h^{\mathrm T}  \left(\mathbf m \times \mathbf m_{xx}\right) -\nu \left(\mathbf m\times \mathbf m_{xx}\right)^{\mathrm T}\left(\mathbf h\times \mathbf m\right)\nonumber\\
& -\nu k  ||\mathbf m\times \mathbf h||_2^2 \label{eqintegrandfirstterm}
\end{align}
and 
%Furthermore, we have that 
\begin{align}
 \mathbf m_{xx}^{\mathrm T} \dot{\mathbf m}   =&-k\mathbf m_{xx}^{\mathrm T} \left( \mathbf m \times \mathbf h\right)+\nu||\mathbf m\times  \mathbf m_{xx} ||_2^2\nonumber\\
  &+\nu k\left(\mathbf m\times  \mathbf h\right)^{\mathrm T} \left(\mathbf m_{xx} \times\mathbf m\right). \label{eqintegrandsecondterm}
\end{align}

Substituting (\ref{eqintegrandfirstterm}) and (\ref{eqintegrandsecondterm}) into equation~(\ref{eqintegrand}) leads to
\begin{align*}
k\mathbf h^{\mathrm T}\dot{\mathbf m} - \mathbf m_{xx}^{\mathrm T} \dot{\mathbf m}
=&2\nu k \left(\mathbf m\times \mathbf m_{xx}\right)^{\mathrm T}\left(\mathbf m\times \mathbf h\right)\\
&-\nu k^2  ||\mathbf m\times \mathbf h||_2^2 -\nu||\mathbf m\times  \mathbf m_{xx} ||_2^2\\
=& -\nu ||\mathbf m\times \mathbf m_{xx} -k \mathbf m\times \mathbf h||_2^2
\end{align*}
Substituting this expression into equation~(\ref{eqLyapunovFunction}) leads to
\begin{align*}
\frac{dV}{dt}  
=&-\nu  ||\mathbf m\times (\mathbf m_{xx} + \mathbf u)||_{\mathcal L_2^3}^2
\leq 0. 
\end{align*}
Thus, the $H_1$-norm of the error does not increase.  \qed
 \end{pf} 

For any equilibrium point $\mathbf r \in E$  of $(\ref{eqcontrolledLLphysical})$ and $\mathbf m \in D$, let $\mathbf m =\mathbf r +\mathbf v$ where $\mathbf v$  is any admissible perturbation; that is $\mathbf v \in D$ and   $\| \mathbf r + \mathbf v\|_2=1. $ 
The linearization of (\ref{eqcontrolledLLphysical}) at $\mathbf r$ is 
\begin{align}
\frac{\partial \mathbf v}{\partial t}&=\nu  \mathbf v_{xx} +\mathbf r \times  \mathbf v_{xx}+ k\mathbf v \times  \mathbf r-k\nu\mathbf r \times\left(\mathbf v \times \mathbf r\right)\nonumber\\
 \mathbf v(0)&=\mathbf v_0.\label{eqcontrolledlinearLLphysical}
\end{align}

\begin{thm}\label{thmlinearr0isasymstablephysical}
Let $\mathbf r \in E $  be an equilibrium point of (\ref{eqcontrolledlinearLLphysical}).  For any positive constant $k$, $\mathbf r$ is a {locally} asymptotically stable equilibrium of (\ref{eqcontrolledlinearLLphysical}) in the $\mathcal L_2^3$-norm. 
\end{thm}

\begin{pf}
For an admissible  $\mathbf v $ {with $||\mathbf v -\mathbf r||_2\leq 2$}, consider the Lyapunov candidate
\begin{equation*}
V(\mathbf v)=\frac{1}{2}\left| \left|\mathbf v-\mathbf r\right|\right|_{\mathcal L_2^3}^2.
\end{equation*}
It is clear that $V\geq 0$ for all $\mathbf v\in D(A)$ and furthermore, $V(\mathbf v)=0$ only when $\mathbf v=\mathbf r$.  Therefore, $V(\mathbf v)>0$ for all $\mathbf v \in D(A)\backslash  \{\mathbf r\}$.

Taking the derivative of $V(\mathbf v)$ implies
\begin{align*}
\frac{dV}{dt}&=\int_0^L(\mathbf v-\mathbf r)^{\mathrm T}\dot{ {\mathbf v}} dx
\end{align*}
and substituting in (\ref{eqcontrolledlinearLLphysical}) leads to
\begin{align*}
\frac{dV}{dt}=&\nu\int_0^L(\mathbf v-\mathbf r)^{\mathrm T} \mathbf v_{xx}  dx + \nu\int_0^L(\mathbf v-\mathbf r)^{\mathrm T}(\mathbf r \times  \mathbf v_{xx}) dx \\
&+k\int_0^L(\mathbf v-\mathbf r)^{\mathrm T} (\mathbf v \times  \mathbf r) dx\\
& -k\nu \int_0^L(\mathbf v-\mathbf r)^{\mathrm T} \left( \mathbf r \times\left(\mathbf v \times \mathbf r\right)\right)dx.
\end{align*}

The second integral is zero by Lemma~\ref{lemmazerointegral}, and applying integration by parts the first term becomes
\[
-\nu\int_0^L\mathbf v_x^{\mathrm T} \mathbf v_{x}  dx.
\] 
It follows that 
\begin{align*}
\frac{dV}{dt}=&-\nu\int_0^L\mathbf v_x^{\mathrm T} \mathbf v_{x}  dx+k\int_0^L(\mathbf v-\mathbf r)^{\mathrm T} (\mathbf v \times  \mathbf r) dx \\
&-k\nu \int_0^L(\mathbf v-\mathbf r)^{\mathrm T} \left( \mathbf r \times\left(\mathbf v \times \mathbf r\right)\right)dx\\
=&-\nu || \mathbf v_x||^2_{\mathcal L_2^3}+k\int_0^L\mathbf h^{\mathrm T} (\mathbf h \times  \mathbf r) dx\\
&-k\nu \int_0^L\mathbf h^{\mathrm T} \left( \mathbf r \times\left(\mathbf h \times \mathbf r\right)\right)dx
\end{align*}
where $\mathbf h=\mathbf v -\mathbf r$. The middle integral is zero since $\mathbf h^{\mathrm T} (\mathbf h \times  \mathbf r)  = \mathbf r^{\mathrm T} (\mathbf h \times  \mathbf h) =0$, and the last integral can be simplified using the fact that 
\[
 \mathbf h^{\mathrm T}  \left( \mathbf r \times\left(\mathbf h \times \mathbf r\right)\right)= \left(\mathbf h \times \mathbf r\right)^{\mathrm T}  \left(\mathbf h \times \mathbf r\right).
  \]
 Therefore,
  \[
  \frac{dV}{dt}=-\nu || \mathbf v_x||^2_{\mathcal L_2^3}-k\nu ||\mathbf h \times \mathbf r||^2_{\mathcal L_2^3}.
  \]
For $k>0$, 
  \[
  \frac{dV}{dt}=-\nu \left(|| \mathbf v_x||^2_{\mathcal L_2^3}+k||\mathbf h \times \mathbf r||^2_{\mathcal L_2^3}\right)\leq 0
  \]
  and furthermore, ${dV}/{dt} =0 $ if and only if $\mathbf v_x=\mathbf 0$ and $\mathbf h \times \mathbf r =\mathbf v \times \mathbf r=\mathbf 0$. This is true only if $\mathbf v =\alpha\mathbf r$ where $\alpha$ is any scalar. Since $\mathbf v$ must be a constant satisfying $\mathbf v \times \mathbf r=\mathbf 0$, then $\mathbf v = \alpha  \mathbf r$ for some constant $\alpha.$  But since $||\mathbf r+\mathbf v||_2=1$, then $\alpha =0$.
  
It follows that  $\mathbf r$ is a locally asymptotically stable equilibrium point of (\ref{eqcontrolledlinearLLphysical}).\qed
\end{pf}

%%%%%%%
%Example
%%%%%%
\section{Example}\label{secExample}
%\achg{most of this section has been changed}

%\chg{rewrote section and captions}

\begin{figure}[h]
\scalebox{0.56}{
\includegraphics{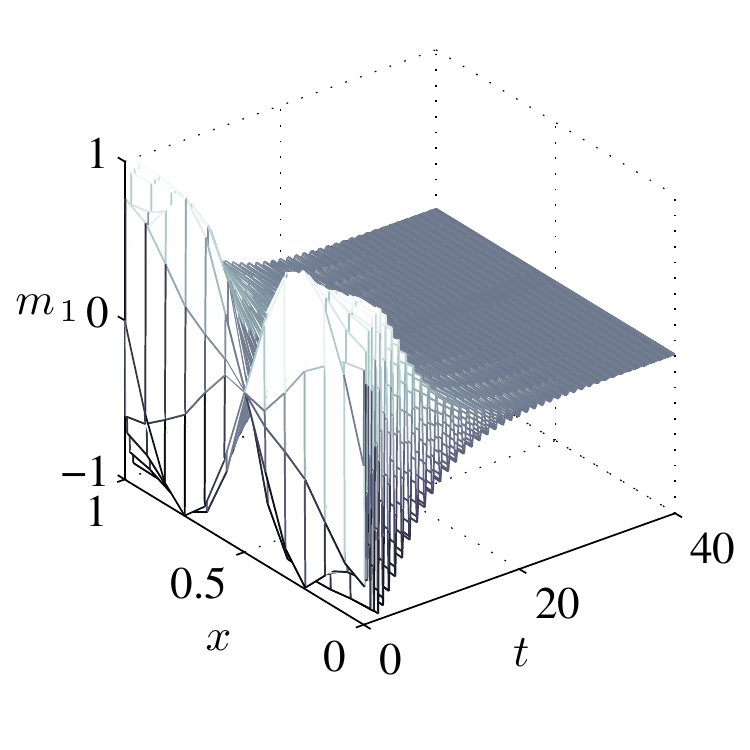}
\includegraphics{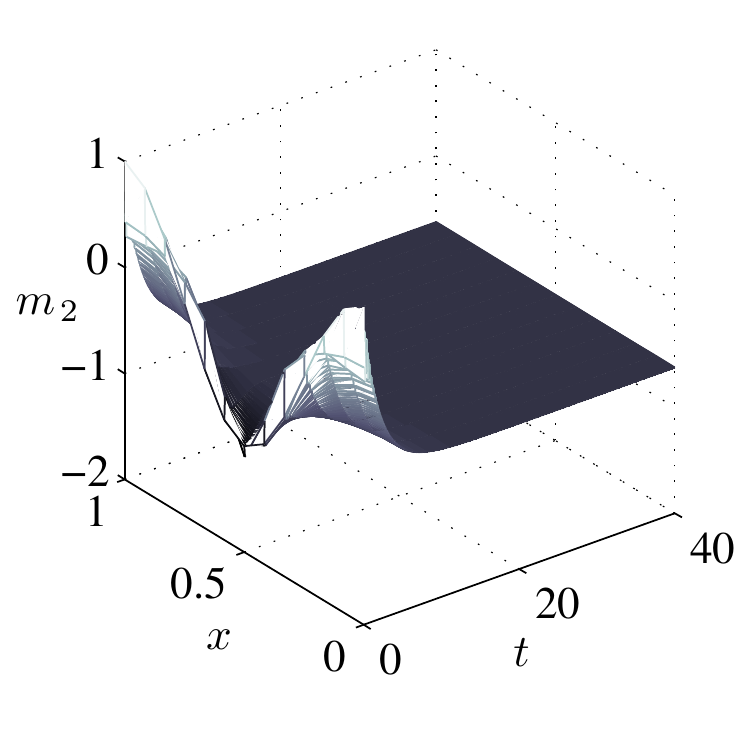}}\\
\scalebox{0.56}{\centering\includegraphics{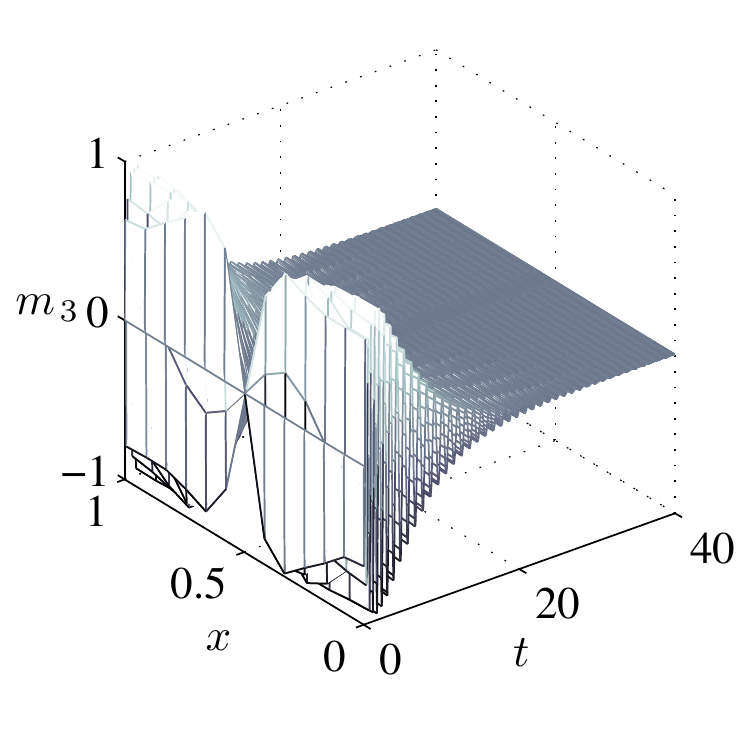}}
 \caption{\label{figMagnetizationNoControl3d}
   Magnetization in the uncontrolled  Landau--Lifshitz equation  moves from initial condition  $\mathbf{{m}_0} (x)$,  to the  equilibrium $\mathbf {r_0}=(0,-0.6,0)$.%   Magnetization in the uncontrolled  Landau--Lifshitz equation  moves from initial condition  $\mathbf{\tilde{m}_0} (x)=(\sin(2\pi x),\cos(2\pi x),0)$,  to the  equilibrium $\mathbf {r_0}=(0,-0.6,0)$.
   }
\end{figure}

Simulations illustrating the stabilization of the Landau-Lifshitz equation were done using  a Galerkin approximation with 12  linear spline elements. For the following simulations, the parameters  are $\nu=0.02$ and $L=1$ with  initial condition  $\mathbf {{m}_0} (x)=(\sin(2\pi x),\cos(2\pi x),0)$. Figure~\ref{figMagnetizationNoControl3d} illustrates that the solution to the uncontrolled Landau--Lifshitz equation settles to $\mathbf{r_0}=(0,-0.6,0).$ 

\begin{figure}[h]
\scalebox{0.56}{
\includegraphics{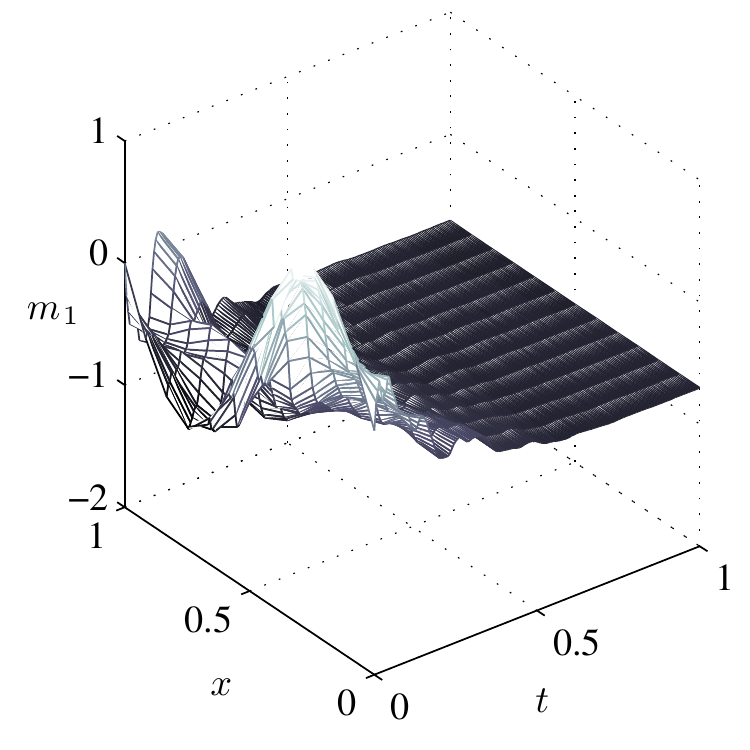}
\includegraphics{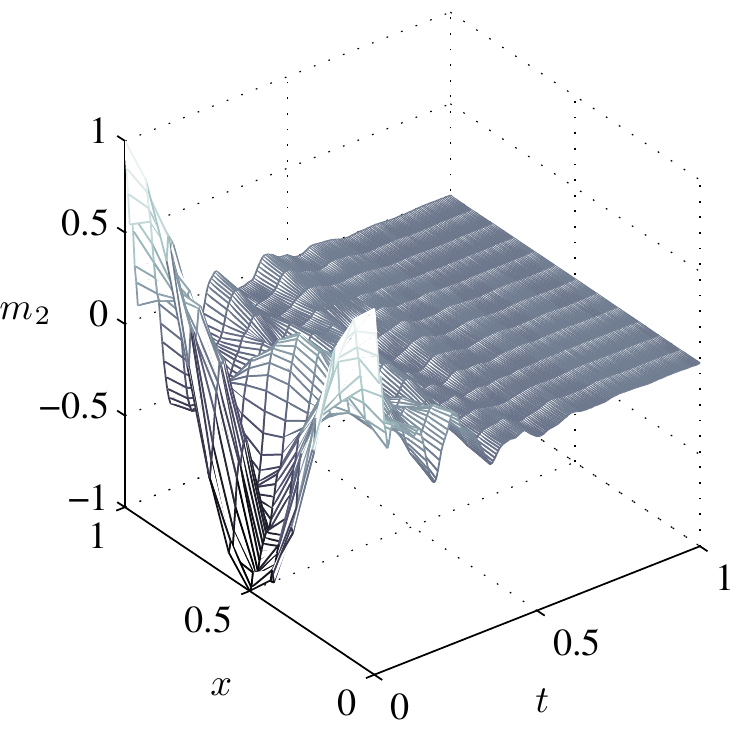}}\\
\scalebox{0.56}{
\includegraphics{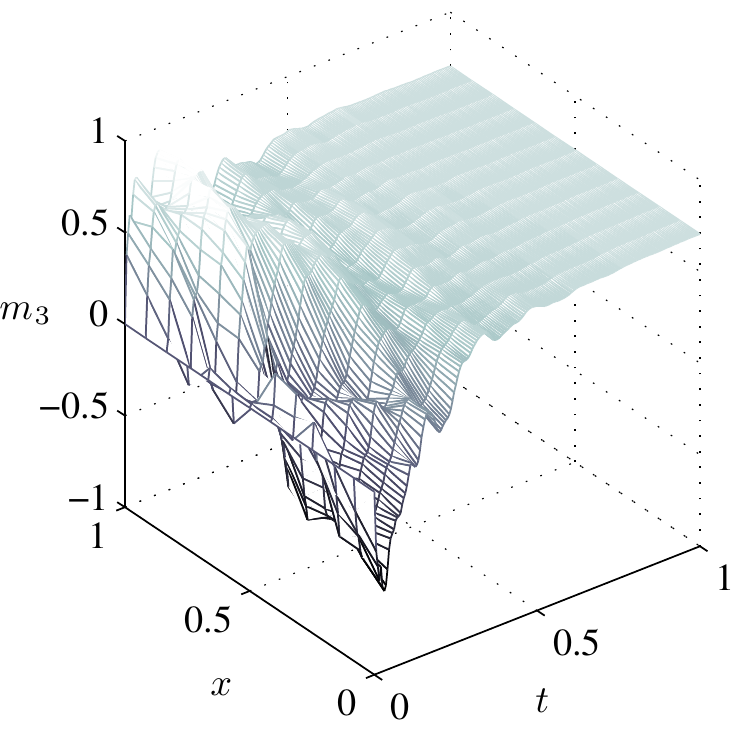}}
 \caption{\label{figMagnetization3d} With a proportional control $(k=0.5)$,  magnetization in the Landau--Lifshitz equation with a linear control moves  from the initial condition
 $\mathbf{{m}_0} (x)$  to the specified equilibrium
 $\mathbf{ r_1}=(-\frac{1}{\sqrt{2}},0,\frac{1}{\sqrt{2}}) .$%With proportional control,  magnetization in the controlled Landau--Lifshitz equation  moves  from the initial condition$\mathbf{\tilde{m}_0} (x)=(\sin(2\pi x),\cos(2\pi x),0)$  to the specified equilibrium $\mathbf{ r_1}=(-\frac{1}{\sqrt{2}},0,\frac{1}{\sqrt{2}}) .$
  }
\end{figure}

\begin{figure}[h]
%\scalebox{0.56}{
%\includegraphics{LLm1controlTWICEICEQEQMeshBone.pdf}
%\includegraphics{LLm2controlTWICEICEQEQMeshBone.pdf}}\\
%\scalebox{0.56}{
%\includegraphics{LLm3controlTWICEICEQEQMeshBone.pdf}}
\scalebox{0.56}{
\includegraphics{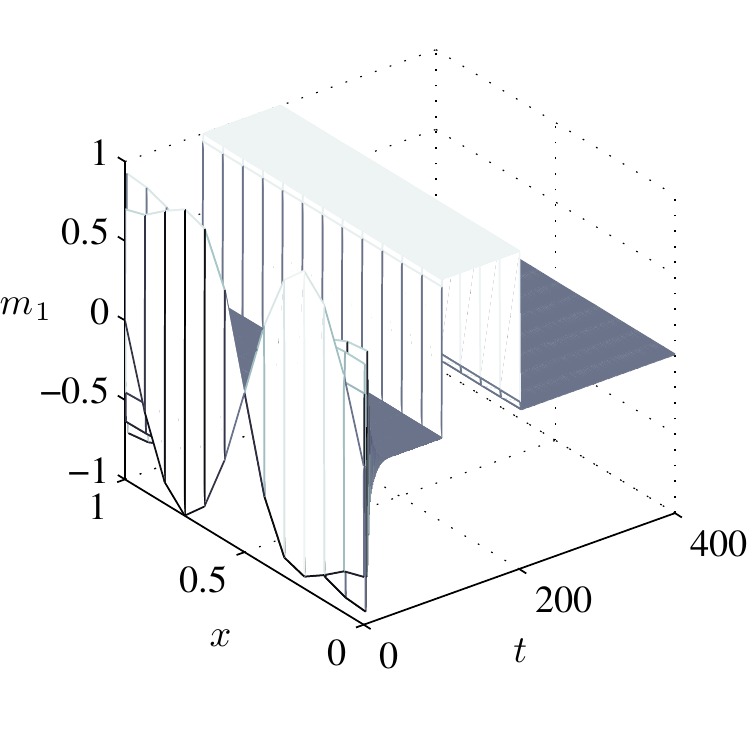}
\includegraphics{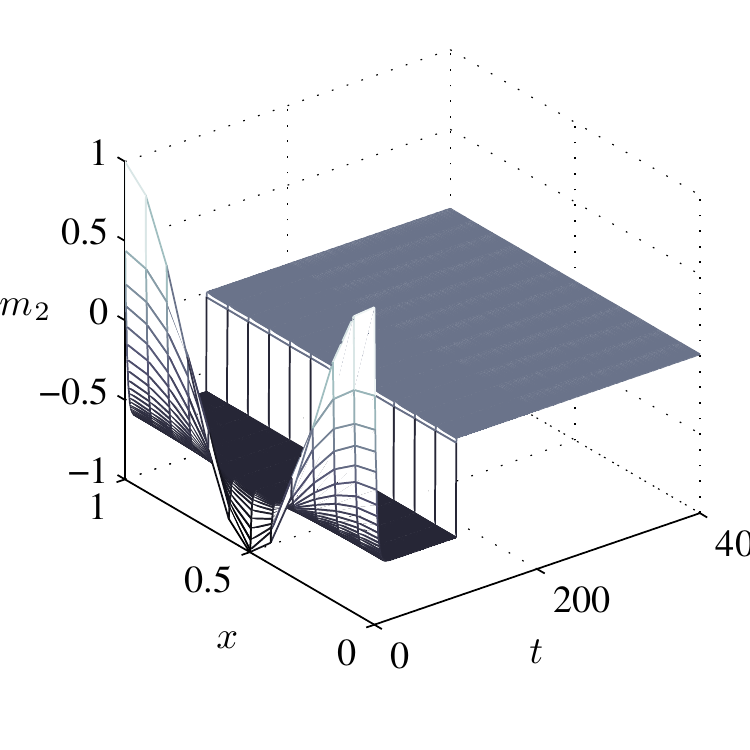}}\\
\scalebox{0.56}{
\includegraphics{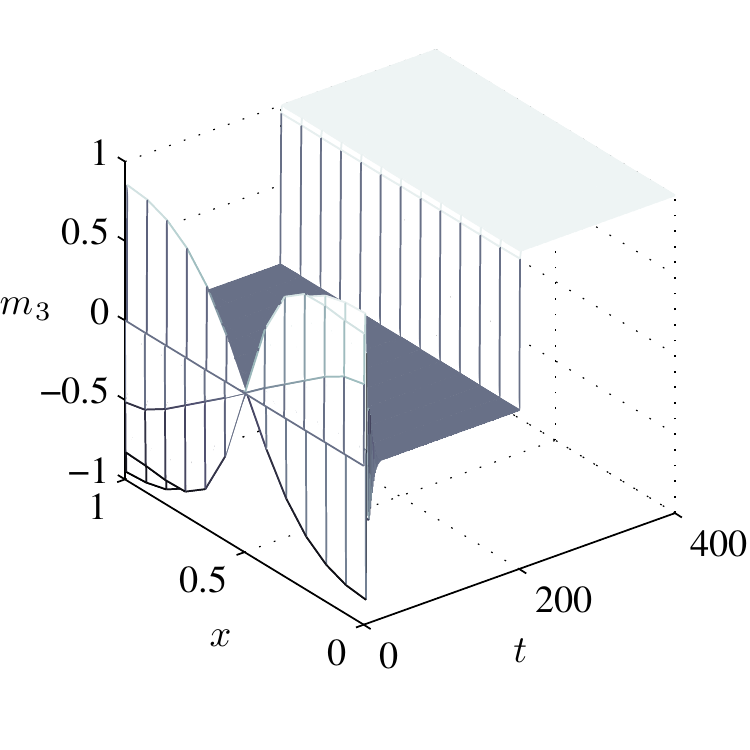}}
 \caption{\label{figLLmcontrolTWICEICEQEQ3d}  Steering magnetization between specified equilibria with a linear control. The uncontrolled magnetization moves from initial condition $\mathbf{{m}_0}$  to $\mathbf {r_0}=(0,-0.6,0)$.   Proportional control  $(k=0.5)$ with two successive values of $\mathbf r$ first forces the magnetization to $\mathbf {r_2}= (1,0,0)$ and then to $\mathbf{r_3}= (0,0,1)$. 
 %($L=1$, $\nu=0.02$, $\mathbf m_0(x)=(\sin(2\pi x),\cos(2\pi x),0)$, $b=1$, $k=0.5$) 
 }
\end{figure}

%ADD Landau--Lifshitz equation + Linear Controlled 
Stabilization of the Landau--Lifshitz equation with affine control (\ref{eqcontrolledLL}) is illustrated in Figure~\ref{figMagnetization3d} with the second equilibrium point chosen to be $\mathbf r_1=(-\frac{1}{\sqrt{2}},0,\frac{1}{\sqrt{2}})$. The control parameter  is $k=0.5$.  Figure~\ref{figLLmcontrolTWICEICEQEQ3d} depicts applying the control twice in succession,   forcing the system from the equilibrium $\mathbf{r_0}$  to  $\mathbf r_2=(1,0,0)$  and then to a new equilibrium $\mathbf r_3 = (0,0,1)$. In each case,  the state of the controlled system  converges to the specified point $\mathbf r_i $ as predicted by the analysis.

\begin{figure}[h]
\scalebox{0.56}{
\includegraphics{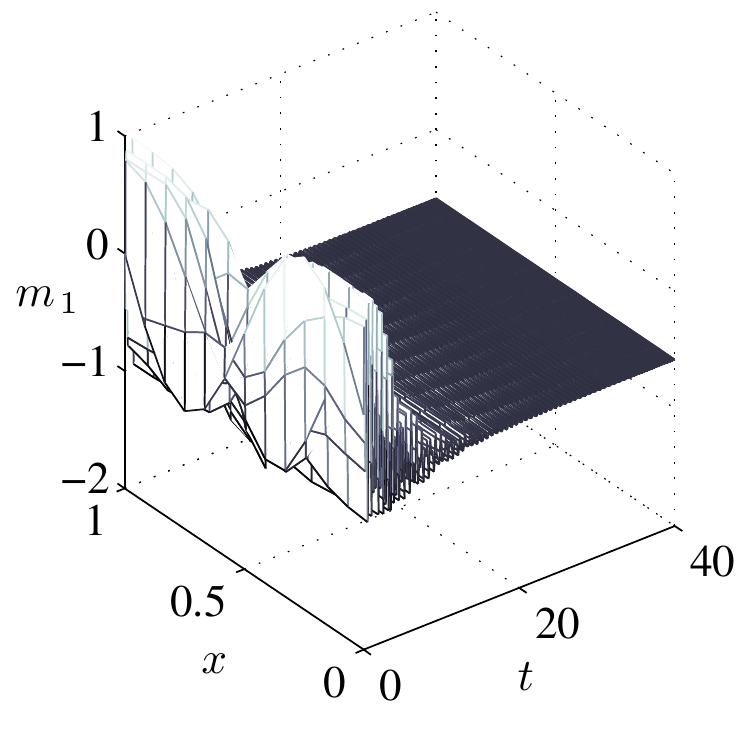}
\includegraphics{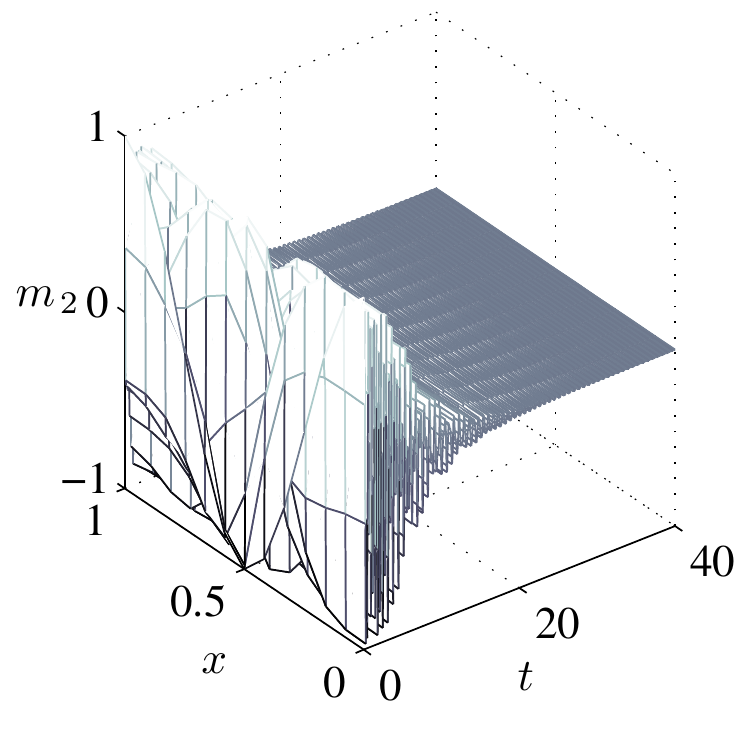}}\\
\scalebox{0.56}{
\includegraphics{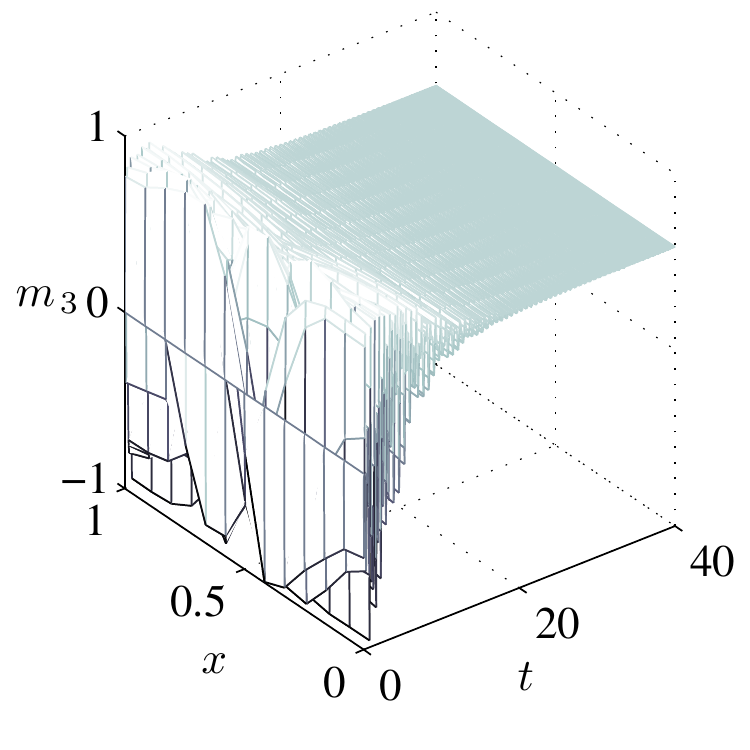}}
 \caption{\label{figMagnetization3dPhysical}  Magnetization in the Landau--Lifshitz equation with nonlinear control moves  from the initial condition   $\mathbf{{m}_0} (x)$  to the specified equilibrium $\mathbf{ r_1}=(-\frac{1}{\sqrt{2}},0,\frac{1}{\sqrt{2}}) $ with control parameter $k=10$. }
\end{figure}

\begin{figure}[h]
\scalebox{0.57}{
\includegraphics{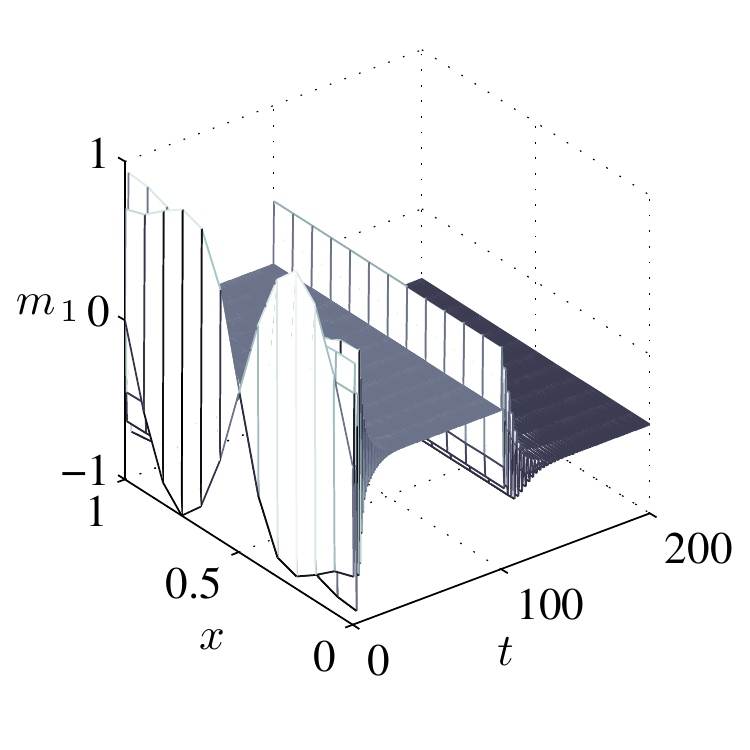}
\includegraphics{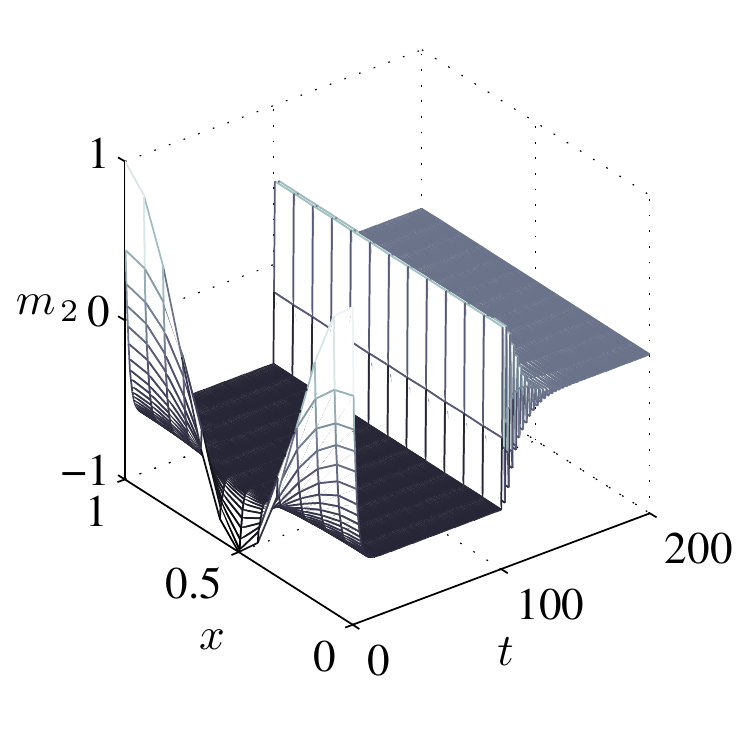}}\\
\scalebox{0.57}{
\includegraphics{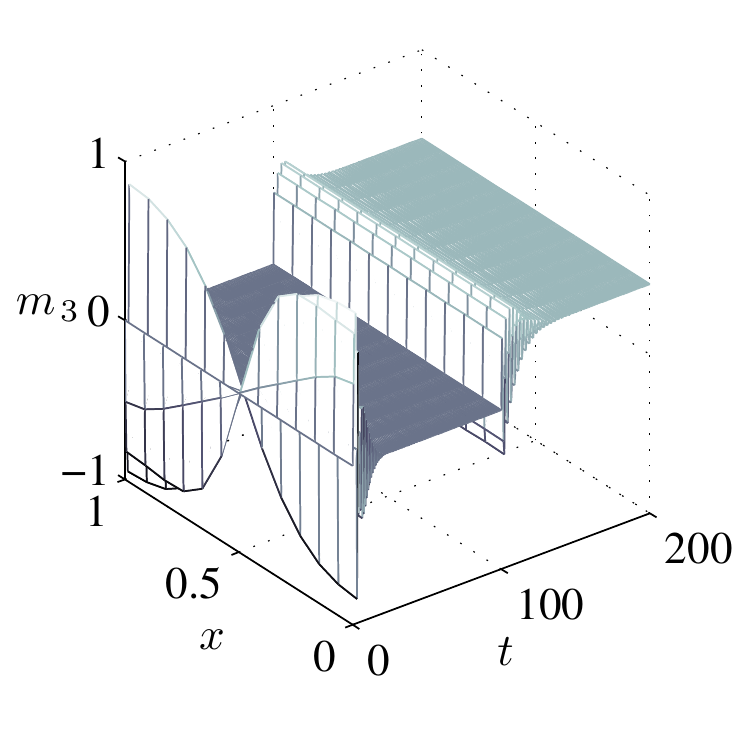}}
 \caption{\label{figMagnetization3dPhysicalICtoEQtoEQ}  Magnetization in the Landau--Lifshitz equation  moves  from the initial condition 
 $\mathbf{{m}_0} (x)$  to the equilibrium $\mathbf {r_0}=(0,-0.6,0)$ without a control. The control, $\mathbf u$ with $k=10$ is then applied to the equation nonlinearly and steers the dynamics to the specified equilibrium
 $\mathbf{ r_1}=(-\frac{1}{\sqrt{2}},0,\frac{1}{\sqrt{2}})$. }
\end{figure}

\begin{figure}[h]
\scalebox{0.57}{
\includegraphics{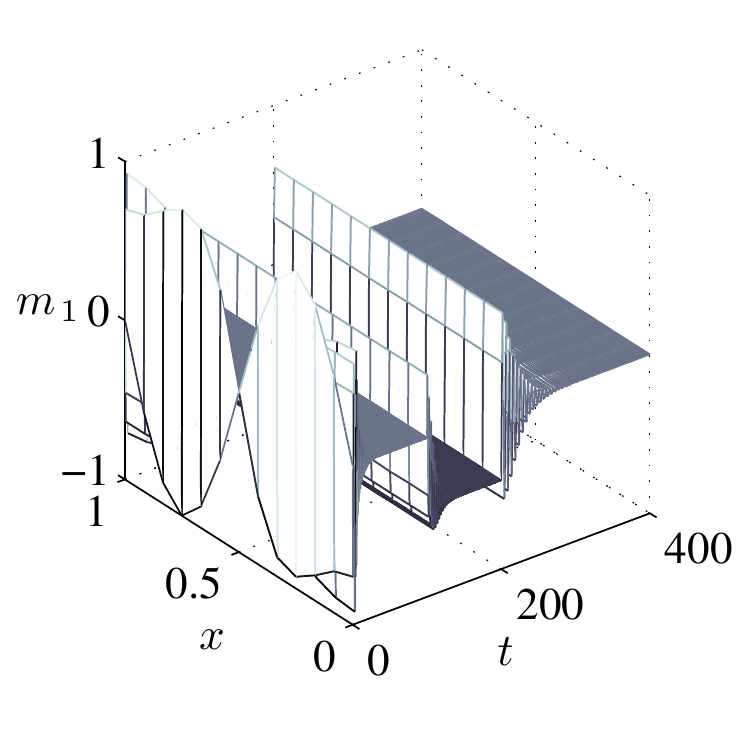}
\includegraphics{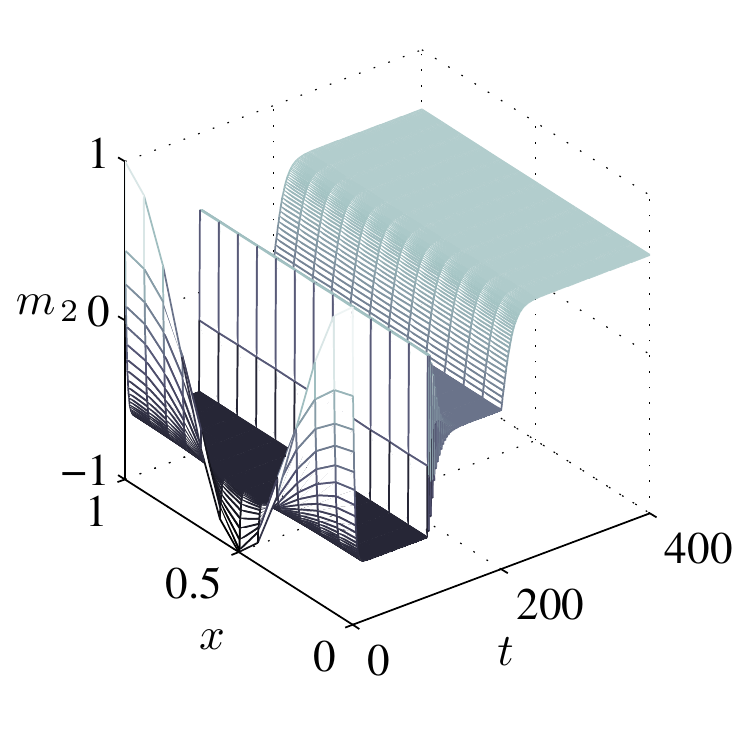}}\\
\scalebox{0.57}{
\includegraphics{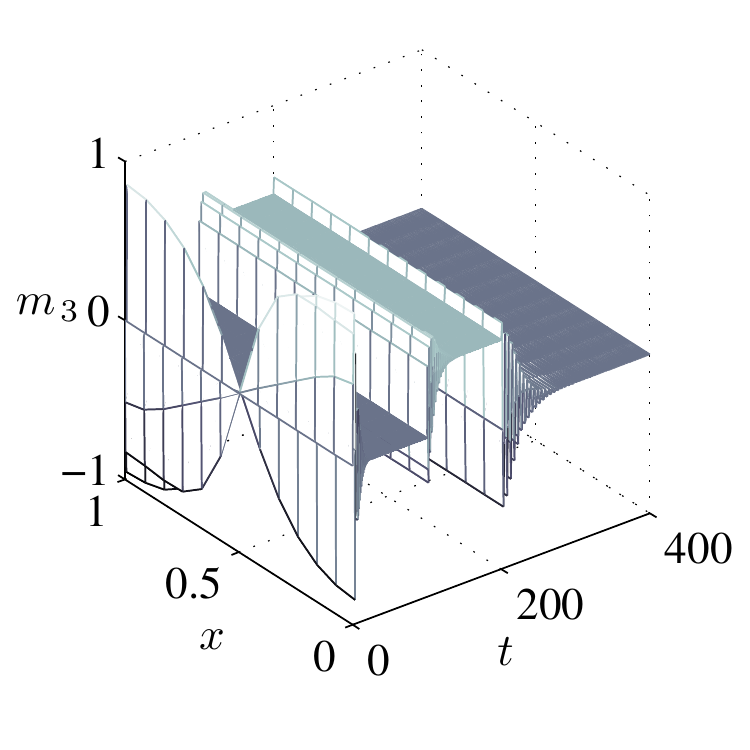}}
 \caption{\label{figMagnetization3dPhysicalICtoEQtoEQtoEQ}  Magnetization in the Landau--Lifshitz equation  moves  from the initial condition 
 $\mathbf{{m}_0} (x)$  to the equilibrium $\mathbf {r_0}=(0,-0.6,0)$ without a control. The control $\mathbf u$ with $k=10$ is then applied to the equation nonlinearly and steers the dynamics to the specified equilibrium
 $\mathbf{ r_1}=(-\frac{1}{\sqrt{2}},0,\frac{1}{\sqrt{2}}) $, and then to another equilibrium, $\mathbf{ r_4}=(0,1,0) $. }
\end{figure}

Stabilization of the Landau--Lifshitz equation with nonlinear control (\ref{eqcontrolledLLphysical}) is illustrated in Figure~\ref{figMagnetization3dPhysical} with the second equilibrium point chosen to be $\mathbf r_1$. The control parameter is $k=10$. It is clear from the figure that the system converges to the specified equilibrium point, $\mathbf r_1$. The control can also be applied after the dynamics have settled to $\mathbf r_0$ as shown in Figure~\ref{figMagnetization3dPhysicalICtoEQtoEQ}. In Figure~\ref{figMagnetization3dPhysicalICtoEQtoEQtoEQ}, the dynamics settle (without a control) to $\mathbf r_0$, then the control is applied in succession twice, which forces the system from $\mathbf r_0$ to $\mathbf r_1$, and then finally to $\mathbf r_4 = (0,1, 0)$. The rate of convergence is slower and the value of $k$ needed is larger than with an affine control. %This is illustrated in Figure~\ref{figMagnetization3dPhysical}.

%%%%%%%%
%Conclusion
%%%%%%%%
\section{Conclusion}
%\achg{most of this section has been changed}

The Landau-Lifshitz equation is a nonlinear system of PDEs with multiple equilibrium points. The fact that it is not quasi-linear means that the linearization  is not guaranteed to predict stability of the non-linear equation \cite{alJamal2013}.  Furthermore, since the objective of the control is to steer between equilibrium points, a linearized analysis, which only yields local results, would not predict stability of the controlled system.  However,  the presence of a 0 eigenvalue in the linearized equation suggested that simple feedback proportional control could  be used to steer the system to an arbitrary equilibrium point. 
 
The  controlled system with an affine control term was shown to be well-posed and also globally asymptotically stable.  Furthermore, it is exponentially stable in the $H_1$-norm. Analysis of the linearized Landau--Lifshitz equation, which shows that the equilibrium is exponentially stable, suggests the  original nonlinear Landau--Lifshitz equation is  locally exponentially stable in the $L_2$-norm.  Future research aims to establish this.  

In applications, the control enters through an applied field and the control enters nonlinearly.  It was shown that the Landau-Lifshitz equation with a nonlinear control has a stable equilibrium point and the linearization has an asymptotically stable equilibrium point.  Simulations indicate that proportional control also stabilize the fully nonlinear model.  This suggests the nonlinear equation has an asymptotically stable equilibrium point. Proving this remains an open research problem. 

%%-----------------------------
%%      your bibliography
%%-----------------------------
\bibliographystyle{ieeetr}
\bibliography{ref}

%no control

%LINEAR CONTROL

%%%%
%ArXiv will not work using LLm1controlTWICEICEQEQMeshBone.pdf, LLm2controlTWICEICEQEQMeshBone.pdf and LLm3controlTWICEICEQEQMeshBone.pdf

%PHYSICAL CONTROL

\end{document}